\documentclass[pagenos]{jktr}

\usepackage{epsfig}
\usepackage{amstex}

\newcommand{\cL}{{\mathcal L}}
\newcommand{\cR}{{\mathcal R}}

\newcommand{\cI}{{\mathcal I}}

\newcommand{\cV}{{\mathcal V}}
\newcommand{\C}{{\hspace{-0.3pt}\mathbb C}\hspace{0.3pt}}
\newcommand{\g}{SL_2(\C)}

\newlength{\widi}
\newlength{\widii}

\newcommand{\rcr}{\raisebox{-5pt}{\mbox{}\hspace{1pt}
                  \epsfig{file=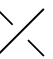}\hspace{1pt}\mbox{}}}
\newcommand{\dbp}{\raisebox{-5pt}{\mbox{}\hspace{1pt}
                  \epsfig{file=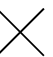}\hspace{1pt}\mbox{}}}
\newcommand{\lcr}{\raisebox{-5pt}{\mbox{}\hspace{1pt}
                  \epsfig{file=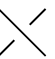}\hspace{1pt}\mbox{}}}
\newcommand{\ift}{\raisebox{-5pt}{\mbox{}\hspace{1pt}
                  \epsfig{file=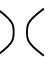}\hspace{1pt}\mbox{}}}
\newcommand{\zer}{\raisebox{-5pt}{\mbox{}\hspace{1pt}
                  \epsfig{file=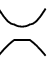}\hspace{1pt}\mbox{}}}

\newcommand{\adb}{\raisebox{-5pt}{\mbox{}\hspace{1pt}
                  \epsfig{file=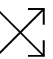}\hspace{1pt}\mbox{}}}
\newcommand{\azr}{\raisebox{-5pt}{\mbox{}\hspace{1pt}
                  \epsfig{file=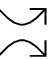}\hspace{1pt}\mbox{}}}
\newcommand{\ain}{\raisebox{-5pt}{\mbox{}\hspace{1pt}
                  \epsfig{file=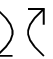}\hspace{1pt}\mbox{}}}
\newcommand{\pls}{%
   \settowidth{\widi}{$sgn(x)=+1$}
   \settowidth{\widii}{$\alpha$}
   \parbox{\widi}{%
         \mbox{}\hfill\raisebox{0pt}[0pt][6pt]{%
                      \makebox[27pt]{$\beta$}}
            \hspace{2pt}\hspace{\widii}\hfill\mbox{}\newline
         \mbox{}\hfill\epsfig{file=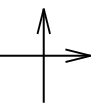}\hspace{2pt}
            \raisebox{11pt}{$\alpha$}\hfill\mbox{}\newline
         $sgn(x)=+1$}}
\newcommand{\mns}{%
   \settowidth{\widi}{$sgn(x)=-1$}
   \settowidth{\widii}{$\alpha$}
   \parbox{\widi}{%
         \mbox{}\hfill\raisebox{0pt}[0pt][6pt]{%
                      \makebox[27pt]{$\beta$}}
            \hspace{2pt}\hspace{\widii}\hfill\mbox{}\newline
         \mbox{}\hfill\epsfig{file=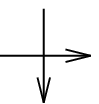}\hspace{2pt}
            \raisebox{11pt}{$\alpha$}\hfill\mbox{}\newline
         $sgn(x)=-1$}}
\newcommand{\xpl}{%
   \settowidth{\widi}{$\alpha$}
   \setlength{\widii}{109pt}
   \addtolength{\widii}{\widi}
   \parbox{\widii}{%
         \mbox{}\hspace{\widi}\hspace{4pt}\raisebox{0pt}[0pt][6pt]{%
             \makebox[27pt]{$\beta$}}\newline
         \mbox{}\hfill\raisebox{11pt}{$\alpha$}\hspace{2pt}
                      \epsfig{file=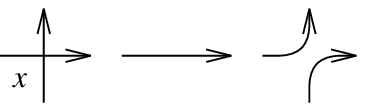}\hfill\mbox{}}}

\newcommand{\rso}{\raisebox{-22pt}[50pt][0pt]{\mbox{}
     \hspace{5pt}\epsfig{file=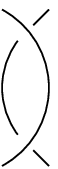}\hspace{5pt}\mbox{}}}
\newcommand{\rsi}{\raisebox{-22pt}[50pt][0pt]{\mbox{}
     \hspace{5pt}\epsfig{file=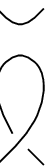}\hspace{5pt}\mbox{}}}
\newcommand{\rsii}{\raisebox{-22pt}[50pt][0pt]{\mbox{}
     \hspace{5pt}\epsfig{file=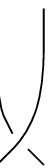}\hspace{5pt}\mbox{}}}
\newcommand{\rsiii}{\raisebox{-22pt}[50pt][0pt]{\mbox{}
     \hspace{5pt}\epsfig{file=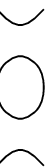}\hspace{5pt}\mbox{}}}
\newcommand{\rsiv}{\raisebox{-22pt}[50pt][0pt]{\mbox{}
     \hspace{5pt}\epsfig{file=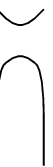}\hspace{5pt}\mbox{}}}
\newcommand{\rsv}{\raisebox{-22pt}[50pt][0pt]{\mbox{}
     \hspace{5pt}\epsfig{file=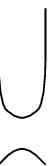}\hspace{5pt}\mbox{}}}
\newcommand{\rsvi}{\raisebox{-22pt}[50pt][0pt]{\mbox{}
     \hspace{5pt}\epsfig{file=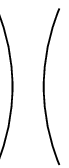}\hspace{5pt}\mbox{}}}

\newtheorem{theorem}{Theorem}

\newtheorem{cor}{Corollary}
\newtheorem{prop}{Proposition}
\newtheorem{conjecture}{Conjecture}
\newtheorem{rem}{Remark}
\newtheorem{sch}{Scholium}

\title{Understanding the Kauffman Bracket Skein Module}
\makeatletter
\author{doug bullock}

\address{Department of Mathematics, Boise State University, Boise, ID
  83725, USA\\ email:  bullock@math.idbsu.edu}

\author{CHARLES FROHMAN}

\address{Department of Mathematics, University of Iowa, Iowa City, IA
  52245, USA\\ email: frohman@math.uiowa.edu}

\author{JOANNA  KANIA-BARTOSZY\'{N}SKA}

\address{Department of Mathematics, Boise State University, Boise, ID
  83725, USA\\ email: kania@math.idbsu.edu}
\makeatother

\begin{document}

\maketitle

\begin{abstract}
The Kauffman bracket skein module $K(M)$ of a 3-manifold $M$ is
defined over formal power series in the variable $h$ by letting
$A=e^{h/4}$. For a compact oriented surface $F$, it is shown that $K(F
\times I)$ is a quantization of the $\g$-characters of the fundamental
group of $F$, corresponding to a geometrically defined Poisson
bracket. Finite type invariants for unoriented knots and links are
defined.  Topologically free Kauffman bracket modules are shown to
generate finite type invariants.  It is shown for compact $M$ that
$K(M)$ can be generated as a module by  cables on a finite set of
knots. Moreover, if $M$ contains no incompressible surfaces, the
module is finitely generated.
\end{abstract}

\section{Introduction}

There are two important classes of invariants of links and manifolds
that are clearly related, however, conceptual bridges connecting the
two are rather sparse.  The first type, classical invariants, are
derivable from the representation theory of $3$-manifold
groups. Examples include the Casson-Walker invariant and the
Chern-Simons invariant.  Quantum invariants are the other type.  Using
a difficult axiomatic computation, Murakami showed that if the quantum
$SU(2)$ invariants are expanded as power series, then the linear terms
recapture the Casson-Walker invariant \cite{Mu}.  Witten derived an
asymptotic formula, using path integrals, relating quantum invariants
to classical invariants. Using tools from classical analysis, Rozansky
proved that Witten's formula is essentially correct when the manifolds
are Seifert Fibered spaces with $3$-singular fibers over the sphere
\cite{Ro}.  What is needed is a routine for passing from one kind of
information to the other.

In a series of papers, Alekseev, Buffenoir, Grosse, Roche and
Schomerus \cite{AGS, BR1,BR2,Bu} have approached this problem using
tools from lattice gauge field theory. It is an easy leap from their
work to a construction of a quantization of the $\g$-characters of a
surface group. However, the construction is so abstract that it is
difficult to get a handle on what this quantization ``looks'' like. It
seems it should play a central role in the transference of data from
the classical setting to the quantum setting.  What we propose here,
is a natural quantization of the $\g$-characters of a surface group.
Since there are already constructions of quantum invariants using the
Kauffman bracket skein module \cite{Bl,G,L}, it should come as no
surprise that we use it here.  The Kauffman bracket skein module of $F
\times I$, where $F$ is a compact orientable surface, is a
quantization of the the ring of $\g$-characters of the fundamental
group of $F$ with respect to a geometrically defined complex linear
Poisson structure.

There is already a small body of literature on skein quantizations of
Poisson algebras \cite{HP,T} associated to surfaces.  
Much of this work focuses on algebras
which lack the refined structure of the ring of characters.  In
particular, one finds quantizations of a symmetric tensor algebra
constructed from loops on the surface.  These use skein modules other than
the Kauffman bracket module, and somewhat different definitions of
quantization.

Przytycki was aware that the Kauffman bracket skein module of $F
\times I$, defined over Laurent polynomials in the variable A can be
viewed as a deformation of the commutative algebra obtained by setting
$A=-1$. Furthermore, he knew that this module is free.  
We show that a completion of it is actually a natural quantization of the
$\g$-characters of a surface group with respect to the standard
Poisson structure coming from a bilinear form on $sl_2$ \cite{BG,Go}.

This paper rests heavily on the  work of \cite{B1,B2,B3,PS} tying
the Kauffman bracket skein module to the ring of $\g$-characters
of the fundamental group of a $3$-manifold.
We exploit the isomorphism  from  \cite{B1} in
order to relate state sum invariants with invariants coming from
representation theory. The concept that makes this paper work
is redefining  the Kauffman bracket
skein module over formal power series.  This is similar in spirit to
Murakami's expansion of $SU(2)$ invariants as  power series.  In our
expansion, the 0-th order term is the ring of $\g$-characters.   
The authors expect that skein
modules will  become a unifying concept in the theory of invariants of 
$3$-manifolds.
  
In the next section Poisson rings, quantization, the Kauffman bracket skein
module, and the ring $X[M]$ of $\g$-characters of a manifold $M$ are
defined. 
In the third section  a geometric  Poisson bracket is defined
 on the ring of  $\g$-characters
of any compact oriented surface $F$.   We prove  that $K(F \times I)$
forms a quantization of $X[F]$ 
with respect to that Poisson bracket. Finally, in the last section we discuss
finite type invariants and topological generators.  In \cite{Pr2}
 Przytycki defined finite
type invariants in a way similar to ours. 

\section{Definitions}

Let $A$ be a commutative algebra with unit over the complex numbers
$\C$. A {\bf Poisson bracket} on $A$ is a bilinear map $ \{\ ,\ \} : A
\otimes A \rightarrow A$ that satisfies three conditions.
\begin{itemize}
\item  The first is antisymmetry, for every  $ a,b \in A$,
$\{a,b\} = - \{b,a\}$.
\item  The second is the Jacobi identity, for every $a,b,c \in A$,
\[ \{a,\{b,c\}\} +\{b,\{c,a\}\} + \{c,\{a,b\}\} = 0 . \]
\item  Finally, the bracket must be
a derivation. Hence, for every $a,b,c \in A$, \[ \{ab,c\} = a\{b,c\}+
 b\{a,c\}. \]
\end{itemize}
 A commutative algebra equipped with a Poisson bracket is called a
{\bf Poisson algebra}.

Denote by $\C[[h]]$ the ring of formal power series in $h$,
topologized by using the sets $a +h^n
\C[[h]]$ as a neighborhood basis for the power series $a$. This is
known as the $h$-adic topology. We work in the category of topological
$\C[[h]]$-modules.  If $V$ is a vector space over $\C$, then we can
form the module $V[[h]]$ of formal power series with coefficients in
$V$ in the obvious way \cite{Ka}. A $\C[[h]]$-module $M$ is {\bf
topologically free} if there exists a vector space $V$ over $\C$ so
that $M$ is isomorphic to $V[[h]]$.

 If $A$ is a Poisson algebra then, following \cite{KS}, we define a
{\bf quantization} of $A$ to be a $\C[[h]]$-algebra $A_h$ together
with a $\C$-algebra isomorphism, $\Theta: A_h/hA_h \rightarrow A$, so
that
\begin{itemize}
\item   as a module $A_h$ is topologically free, and
\item  if 
$a,b \in A$ and $a',b'$ are any elements of $A_h$ with $\Theta(a')=a$
and $\Theta(b')=b$, then \[\Theta\left( \frac{a'b'-b'a'}{h}\right) =
\{a,b\}, \] where the equation above implicitly uses the fact that
$A_h$ is topologically free.
\end{itemize}

\begin{figure}
    \mbox{}\hfill\epsfig{figure=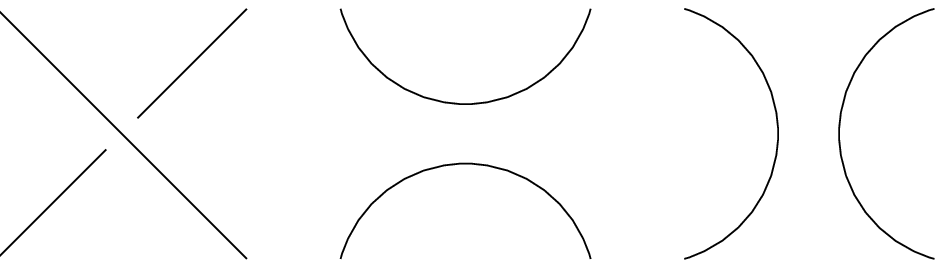,width=2in}\hfill\mbox{}
    \newline\mbox{}\hspace{1.7in}\makebox[1.65in]{$L$\hfill$L_0$\hfill
       $L_\infty$}
    \caption{}\label{relations}
  \end{figure}

The notion of a skein module was first introduced by Przytycki in
\cite{Pr}.  Let $M$ be a compact oriented 3-manifold.  The {\bf
Kauffman bracket skein module} of $M$ is an algebraic invariant, denoted
$K(M)$, which is built from the set $\cL_M$ of framed links in $M$. By
a framed link we mean an embedded collection of annuli considered up
to isotopy in $M$, and we include the empty collection $\emptyset$.
Three links $L$, $L_0$ and $L_\infty$ are said to be {\bf Kauffman bracket
skein related} if they can be embedded identically except in a ball
where they appear as shown in Fig.\
\ref{relations}
(with the blackboard framing).
The notation $L \amalg \bigcirc$ indicates the union of $L$ with an
unlinked, 0-framed unknot. 
      
Let $\C\cL_M$ be the vector space whose basis is $\cL_M$.  As above,
$\C\cL_M[[h]]$ is the $\C[[h]]$-module of power series with
coefficients in $\C\cL_M$.  Let $t$ be the formal power series
$e^{h/4}$, and let $s$ be the series $t^2$.  If $L$, $L_0$ and
$L_\infty$ are Kauffman bracket skein related then $L+tL_{0}+t^{-1}L_{\infty}$
is called a {\bf skein relation}.  For any $L$ in ${\mathcal L}_M$ the
expression $L \amalg \bigcirc + (s+s^{-1})L$ is called a {\bf framing
relation}.  Let $S(M)$ be the smallest closed submodule of
$\C\cL_M[[h]]$ containing all possible skein and framing relations.
We define $K(M)$ to be the quotient $\C\cL_M[[h]]/S(M)$. 

$K_0(M)$ is the quotient $K(M)/hK(M)$.  Two links represent the same
element in $K_0(M)$ if they are homotopic as maps from a disjoint
union of circles into $M$.  Furthermore, $K_0(M)$ is an algebra over
$\C$.  To define multiplication we begin with elements of $K_0(M)$
that can be represented by links, and extend by requiring that
multiplication distribute over addition.  To multiply the equivalence
classes of two links, first perturb them to be disjoint, and then take
the class of their union.

Let $M$ be a compact oriented manifold, with or without boundary. The
set of $\g$-representations of the fundamental group of $M$ can be
viewed as an affine algebraic set $\mathcal{R}(M)$.  Choose a finite
set of generators $\{\gamma_i\}_{i=1}^n$ and relations $r_j$ for the
fundamental group of $M$. The set of points in $\prod_{i=1}^n\g$
satisfying the matrix equations given by the $r_j$ is in one to one
correspondence with the set of representations of $\pi_1(M)$ into
$\g$. The sets obtained from different presentations of $\pi_1(M)$ are
algebraically equivalent \cite{CS}.  Let $R[M] $ be the coordinate
ring of $\mathcal{R}(M)$. Conjugation of representations induces an
action on the ring $R[M]$.  If $f \in R[M]$, $A \in \g$ and $\rho :
\pi_1(M) \rightarrow \g$ then we define $A.f$ to be the function with
$A.f(\rho)= f(A^{-1}\rho A) $.  The {\bf ring of characters} of $M$,
$X[M]$, is the subring of $R[M]$ fixed by this action. The ring $X[M]$
is the coordinate ring of the $\g$-character set, $\mathcal{X}(M)$.
 
For each $\gamma \in \pi_1(M)$ define a function $tr_\gamma : \cR(M)
\rightarrow \C$ by $tr_{\gamma}(\rho)=tr(\rho(\gamma))$, where $tr$  denotes the standard
trace on  $\g$.
The algebra $X[M]$ is generated by the
functions $tr_{\gamma}$.
Trace on $\g$ satisfies three important identities. 
\begin{itemize}
\item For any $A,  B \in \g$, $tr(AB) =tr(BA)$.
\item For any $A \in \g$, $tr(A) = tr(A^{-1})$.
\item For any  $A, B  \in \g$, $tr(AB) + tr(AB^{-1})=tr(A)tr(B)$.
\end{itemize}

The first property implies that the functions $tr_{\gamma}$ are in the
fixed subring of $R[M]$. It also implies that if $\gamma$ and $\delta$
are conjugate elements of $\pi_1(M)$ then $tr_{\gamma} = tr_{\delta}$.
Since the conjugacy classes of the fundamental group of any space can
be identified with free homotopy classes of maps from the circle into
that space, we can define $tr_{\gamma}$ even when $\gamma$ is a loop
that is not based at the basepoint $x_0$ of the fundamental group. To
do this, use a path from a point on $\gamma$ back to $x_0$ to
conjugate $\gamma$ to be based at $x_0$, then apply the trace. The
second property implies that we don't need to be working with oriented
loops; no matter how you orient a loop $\gamma$ you get the same
function $tr_{\gamma}$.  The third property allows us to define an
isomorphism from $ K_0(M) $ to $X[M]$.  If $\Gamma$ is a union of
loops $\gamma_i$, $i=1 \ldots n$, then we let $n_{\Gamma}=
\prod_i-tr_{\gamma_i}$. We can define a map
\[\Theta : K_0(M) \rightarrow X[M]\]
by sending $\Gamma$ to $n_{\Gamma}$ and extending  linearly.    The
map is a well defined surjection of algebras \cite{B1} whose kernel
is the ideal of nilpotent elements of $K_0(M)$. 

We restrict our attention to surfaces, where Przytycki and Sikora
\cite{PS} have shown that there are no nilpotents.  If two loops from
a collection $\Gamma$ intersect transversely, then we can resolve the
intersection via
\begin{equation}\label{charsk}
  \adb = - \azr - \ain
\end{equation}
If $\Gamma_0$ is the first smoothing and
$\Gamma_\infty$ is the second, then (\ref{charsk}) implies 
\begin{equation}
n_{\Gamma}= -n_{\Gamma_0}-n_{\Gamma_{\infty}}
\end{equation}
in $X[M]$.  This allows us to give a state sum model for rewriting
polynomials in the traces of curves in terms of polynomials in the
traces of curves without intersections. Suppose that there are $k$
crossings in $\Gamma$ that you want to remove. There are two
smoothings at each crossing, hence, $2^k$ different ways of smoothing
all the crossings.  We call each of these a {\bf state}. It is easy to
see that
\begin{equation}
 n_{\Gamma} =(-1)^k \sum_{\text{states}\ S} n_S.
\end{equation}

\section{The Kauffman Bracket Module as a Quantization}

We begin with some reminders that will help the reader understand where
the Poisson bracket on the $\g$-characters of a surface group comes from.
The standard basis for the Lie algebra $sl_2$ is given by
\[ X = \begin{pmatrix} 0 & 1  \\
                       0 & 0  \\
       \end{pmatrix}, \qquad
   H = \begin{pmatrix} 1 & 0  \\
                       0 & -1 \\
       \end{pmatrix} \quad \text{and} \quad
   Y = \begin{pmatrix} 0 & 0  \\
                       1 & 0  \\
       \end{pmatrix}.\]
 We use the nondegenerate
$ad$-invariant bilinear form 
 $ B : sl_2 \otimes sl_2 \rightarrow \C$,  which, with respect to the standard
basis, is given by 
\[ \begin{pmatrix} 0 & 0 & -1 \\
0 & -2 & 0 \\
-1& 0 & 0 \end{pmatrix} .\]  It should be noted that $B$ is $\frac{1}{4}$ of
the standard Killing form.
The form $B$ allows us to identify $sl_2^*$ and $sl_2$. Hence, the
complex differential of a holomorphic map $f: \g \rightarrow \C $ has
a {\bf gradient}, $\delta_f$, with respect to $B$.  The gradient is
defined by $\delta_f(A)=Z$, where $B(Z,W)= df_A(W)$ and  the complex  tangent
space of $\g$ at $A$ is identified  with $sl_2$ via left translation.  

The gradient of the trace is an interesting example.  We begin by
computing its differential at the point
\[A=  \begin{pmatrix} a_{11} & a_{12} \\ a_{21} & a_{22} \end{pmatrix}.\]
To do this, choose a path $ \gamma(h)=A \exp(hX)$ through $A$ in the
direction of $X$.  Then
\begin{equation}
dtr_A(X)=\left.\frac{d}{dh}tr(\gamma(h))\right|_{h=0}=a_{21}.
\end{equation}
The values of $dtr_A$ on $Y$ and $H$ are computed by differentiating
along the paths $A \exp(hY)$ and $A \exp(hH)$ respectively.  From the
above definition of gradient we obtain
\begin{equation}\label{grad}
 \delta_{tr}(A) = -a_{21} Y -a_{12}X -\frac{1}{2}(a_{11} - a_{22})H.
\end{equation}
It is also worth noting that
\begin{equation}\label{god}
 B(\delta_{tr}(A),\delta_{tr}(B)) = \frac{1}{2}
 tr(AB^{-1})-\frac{1}{2} tr(AB).
\end{equation}

In \cite{BG}, building on work of \cite{Go}, it is shown that the
coordinate ring of the character variety $\mathcal{X}_G(F)$ of the
representations of the fundamental group of a compact oriented surface
$F$ into a Lie group $G$ admits a Poisson structure.  Let $B$ be an
$ad$-invariant bilinear form on the Lie algebra ${\mathfrak g}$ of
$G$.  We use $B$ to define $\delta_f : G \rightarrow {\mathfrak g}$
for any $f : G \rightarrow {\mathbb R}$ as above. Let $\gamma$ and
$\eta$ be conjugacy classes of $\pi_1(F)$.  Let $P$ and $Q$ be
transverse loops representing $\gamma$ and $\eta$.  If $x \in P \cap
Q$, define the sign of $x$ via the convention shown in Fig.\
\ref{signs}.
\begin{figure}
\mbox{}\hfill\pls\hfill\mns\hfill\mbox{}
\caption{}\label{signs}
\end{figure}
For each point $x \in P \cap Q$, choose a path $c_x$ joining the
basepoint of $\pi_1(X)$ to $x$.  Let $\gamma_x$ be the element of
$\pi_1(X)$ obtained by first traversing $c_x$, then $\gamma$ and then
$c_x^{-1}$.  Define $\eta_x$ similarly.

As with $\g$, we may form the set $\cR_G(F)$ of representations into
$G$.  Then $X_G[F]$ is, again, the fixed subring of its coordinate ring.
If $f: G \rightarrow {\mathbb R} $ is an invariant function on $G$
then we can define $f_{\gamma} \in X[F]$ by letting $f_{\gamma}(\rho) = f(\rho(\gamma))$.    According to
\cite{BG}, there is a Poisson bracket given by the formula
\begin{equation}\label{fish}
\{f_{\gamma}, f'_{\eta}\}(\rho) = \sum_{x \in P \cap Q}
sgn(x) B(\delta_{f}(\rho(\gamma_x)),\delta_{f'}(\rho(\eta_x))).
\end{equation}

We wish to interpret (\ref{fish}) in our own setting. Making the minor
adjustment of working complex linearly, we use the previously defined
form $B$ on $sl_2$. In this case the only invariant function needed is
the trace, whose gradient is given in (\ref{grad}).  For convenience, we
allow one symbol to denote both a curve on a surface and its free
homotopy class.  If $\alpha$ and $\beta$ intersect transversely at a
point $x$, define $\alpha \beta _x$ to be the loop constructed by
first traversing $\alpha$ beginning at $x$ and then traversing $\beta$
(Fig.\
\ref{splice}).
\begin{figure}[b]
  \mbox{}\hfill\xpl\hfill\mbox{}
\caption{}\label{splice}
\end{figure}
Using (\ref{god}) we see that
\begin{equation}\label{wellandgood}
\{ tr_{\alpha},tr_{\beta}\} = \sum_{x \in \alpha \cap \beta} sgn(x)
\frac{1}{2}(tr_{\alpha\beta^{-1}_x} - tr_{\alpha \beta_x}).
\end{equation}

\begin{figure} 
  \begin{center}
  \begin{tabular}{rlcr}
  & \makebox[.7in]{}     & \makebox[.7in]{} & \makebox[.7in]{} \\
   \rule[-10pt]{0pt}{15pt} &\hspace{11pt}$\beta$  &   & \\
  \raisebox{12pt}{$\alpha$\llap} & \epsfig{file=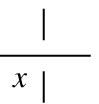} & 
     \epsfig{file=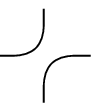} & \epsfig{file=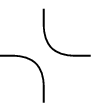} \\
   \rule{0pt}{15pt} &  & $\alpha\beta_{x,0}$ & $\alpha\beta_{x,\infty}$
  \end{tabular}
  \end{center}
\caption{}\label{kauf}
\end{figure}

The formula above is well and good, but it could be better. Recalling that
the trace of $A$ and the trace of $A^{-1}$ are the same, there should be
a formula for the Poisson bracket that does not depend on orientations.
Suppose that
$\alpha$ and $\beta$ are simple closed curves that intersect one another transversely. 
It is simpler to use $n_{\gamma}=-tr_{\gamma}$, because
 the skein rule from the previous
section can be used.
 To represent $\{n_{\alpha},n_{\beta}\}$
draw a link $L_{\alpha \cup \beta}$ consisting of $\alpha$ and $\beta$  with 
$\alpha$
always crossing over $\beta$. At any crossing define two new links,
$\alpha \beta_{x,0}$ and $\alpha \beta_{x,\infty}$ according to the rule
in Fig.\ \ref{kauf}, where you reinterpret the resulting objects as a loop in the
surface.

\begin{prop}
\begin{equation}\label{stillnotgoodenough}
 \{ n_{\alpha},n_{\beta}\}=
 \sum_{x \in \alpha \cap \beta} \frac{1}{2}n_{\alpha \beta_{x,0}}
-\frac{1}{2}n_{\alpha \beta_{x,\infty}}.
\end{equation}
\end{prop}

\begin{proof} 
Compare (\ref{wellandgood}) and (\ref{stillnotgoodenough}).
\end{proof}

The formula above is still not good enough. We need to have a formula for the
Poisson bracket that allows us to compare it to other brackets.
Suppose that $\alpha$ and $\beta$
intersect in $k$ points. 
Each term of the summation for $\{n_{\alpha},n_{\beta}\}$ has $k-1$  crossings
of $\alpha$ over $\beta$ remaining.  There is a state sum model for
(3.9) obtained by resolving all these crossings according to Fig. \ref{kauf}.  
Notice that each
state of this model occurs exactly $k$ times in the resolution. 
Let $0(S)$ be the number of resolutions of type
zero used in making the state $S$, and let $\infty(S)$ be the
number of  resolutions of  type infinity used.

\begin{prop}
\begin{equation}\label{atlast}  
\{ n_{\alpha},n_{\beta}\}=(-1)^k \sum_{\text{states}\ S} \frac{1}{2}
(\infty(S) -0(S))n_S.
\end{equation}
\end{prop}

\begin{proof}
Keep track of the crossings as you resolve.
\end{proof}

Define a Poisson structure on $K_0(F \times I)$ as follows. If $a$ and
$b$ are in $K_0(F \times I)$, choose $\overline{a}$ and $\overline{b}$
in  $K(F \times I)$ that map down to $a$ and $b$ under the quotient map.
Let
\begin{equation}
 \{ a , b \} = \frac{\overline{a} \overline{b} -\overline{b} \overline{a}}
{h}\ \mod \  hK(F \times  I).
\end{equation}

\begin{prop} The isomorphism $\Theta : K_0(F \times I) \rightarrow 
X[F ]$
is a morphism of Poisson algebras. \end{prop}

\begin{proof} Compute $\overline{a} \overline{b} -\overline{b} \overline{a}$
by resolving into states corresponding to the crossings of
$\overline{a}$ and $\overline{b}$. Since the crossings in the first
product are exactly the crossings in the second product we get,
\[\sum_{\text{states}\ S}\left( (-t)^{0(S)-\infty(S)}-(-t)^{\infty(S)-0(S)}\right) S.\]
Notice that the coefficient of $h$ is
\[\frac{1}{2}(-1)^k \sum_{\text{states}\ S} (0(S)-\infty(S)) S,\]
where $k$ is the number of times $\overline{a}$ crosses with
$\overline{b}$. 
\end{proof}

\begin{figure} 
   \mbox{}\hfill
   \epsfig{file=reid0.eps}\hspace{7pt}\raisebox{22pt}{$=$}\hspace{5pt}
   \epsfig{file=reido.eps}\hspace{7pt}\raisebox{22pt}{$=$}\hspace{5pt}
   \epsfig{file=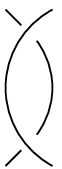}\hfill
   \epsfig{file=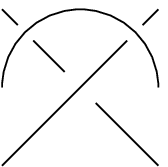}\hspace{7pt}\raisebox{22pt}{$=$}\hspace{5pt}
   \epsfig{file=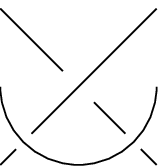}\hfill\mbox{}
\caption{}\label{reide}
\end{figure}

If $K(F \times I)$ is topologically free, then $K(F \times I)$ and
$\Theta$ form a quantization of $X[F]$. The proof that $K(F \times
I)$ is topologically free is an adaptation of ideas of Przytycki
\cite{Pr}.  A diagram of a framed link in $F \times I$ is a possibly
empty four valent graph in $F$, with overcrossings, undercrossings and
blackboard framing. Every framed link can be represented by a diagram,
and every diagram represents a framed link. The correspondence is not
one-to-one.  Two diagrams represent the same framed link if and only if
they differ by a sequence of isotopies of $F$ and Reidemeister moves
of the second and third types (Fig.\ \ref{reide}).
Let ${\mathcal D}(F)$ be the module of power series with coefficients
in the vector space over the set of diagrams up to isotopy.  We define
$RS(F)$ to be the smallest closed submodule of ${\mathcal D}(F)$
containing all skein and framing relations along with the two
Reidemeister moves admissible for regular isotopy.  It is clear that
$K(F \times I)$ is the quotient of ${\mathcal D}(F)$ by $RS(F)$.

A trivial circle on a surface is a circle that bounds a disk.  Let
${\mathcal B}$ be the vector space over $\C$ whose basis consists of
all diagrams that have no crossings and no trivial circles.  Let
${\mathcal B}[[h]]$ be the module of formal power series with
coefficients in ${\mathcal B}$.  Beginning with $\Lambda_0 \in
{\mathcal D}(F)$, construct a sequence as follows.  Given
$\Lambda_{n-1}$, let $\Lambda_{n-1,n-1}$ be the coefficient of
$h^{n-1}$ in $\Lambda_{n-1}$. Resolve the crossings in
$\Lambda_{n-1,n-1}$ via skein relations to get an element $B_{n-1}$ of
${\mathcal B}[[h]]$.  Let $\Lambda_n$ be $\Lambda_{n-1}
-\Lambda_{n-1,n-1} + B_{n-1}$. Notice that, by construction, all of
the $\Lambda_n$ are equivalent modulo $RS(F)$. Furthermore, the
sequence is Cauchy and converges to an element of ${\mathcal B}[[h]]$.
Define the map $\Psi: {\mathcal D}(F) \rightarrow {\mathcal B}[[h]]$
by letting $\Psi(\Lambda )$ be the limit of this sequence when
$\Lambda=\Lambda_0$.  It is clear that $\Psi$ is well defined and
continuous.  Hence if $\Psi(RS(F)) = \{0\}$ then $\Psi$ descends to a
continuous map
\[\Psi : K(F \times I) \rightarrow {\mathcal B}[[h]].\]
Since $\Psi$ is linear and continuous, we only need to check this for
each of the four relations on a single diagram.  The framing relation
is trivial.  We expand the effect of the second Reidemeister move in
Fig.\ \ref{2reid}.  The other two relations work similarly.
\begin{figure}
\begin{equation*}
\begin{split}
\rso & = -t \rsi -t^{-1} \rsii \\
     & = \rsiii +s \rsiv +s^{-1} \rsv + \rsvi \\ 
     & = \rsvi\\
\end{split}
\end{equation*}
\caption{}\label{2reid}
\end{figure}
\par
From these computations we see that there is a continuous map $\Psi :
K(F \times I) \rightarrow {\mathcal B}[[h]]$. There is an obvious map
back obtained by including ${\mathcal B}[[h]]$ into ${\mathcal D}(F)$
and then projecting to the quotient. It is easy to check that these
maps are the inverses of one another. Therefore the two modules are
isomorphic, and we have proved the following theorem.

\begin{theorem}  The Kauffman bracket skein module $K(F \times I)$ and
the map $\Theta : K_0(F \times I) \rightarrow X[F]$ form a quantization
of $X[F]$. \end{theorem}

\section{Applications}

\subsection{Finite type invariants}

One may think of $K(M)$ as a source of invariants of framed,
unoriented links in $M$.  The invariants are most interesting when the module
is topologically free.  Suppose there is an isomorphism between $K(M)$ and
some $V[[h]]$. Given a framed link $L$ in $M$, take the power series
$\sum v_ih^i \in V[[h]]$ corresponding to the class of $L$ in $K(M)$.
 The coefficients $v_i$ are  framed link invariants.
By choosing a basis for $V$ we obtain numerical invariants, namely the
coefficients of each $v_i$ in this basis.

We have seen that $K(F\times I)$ is isomorphic to $V[[h]]$, where $V$
has a natural basis.  When $F$ is a disk these invariants are familiar
ones.  If $L$ is given an orientation then it has a Jones polynomial,
$J_L$, and a well defined writhe, $w(L)$.  Regardless of orientation,
$\langle L \rangle = t^{-3w(L)}J_L(e^h)$, provided $J$ is normalized
at $J_\emptyset = 1$.  It is well known that the coefficients of
$J_L(e^h)$ are (oriented) finite type invariants, and that those of
$t^{-3w(L)}$ are (framed and oriented) finite type invariants.  It
follows that, in some sense, the coefficients of $\langle L \rangle$
are finite type invariants.  We introduce one possible definition that
makes this precise.

\begin{figure}
 \mbox{}\hfill\epsfig{figure=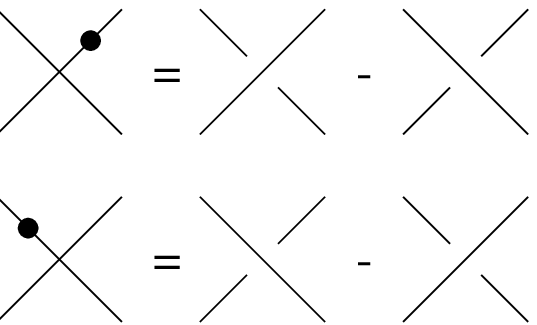,width=2.5in}\hfill\mbox{}
 \caption{}\label{notation} \end{figure}

Recursively, a {\bf decorated link with $n$ double points} represents
a linear combination in $\C\cL_M$. It is a difference of two decorated
links with $n-1$ double points as in Fig.\ \ref{notation}.  Here the
links are assumed to be identical outside the neighborhood shown, and
the arcs have a blackboard framing.  The dot denotes which of the two
possible differences is meant.  A choice of dots for the double points
of a link is called a {\bf decoration}.
 
Noting that a different decoration can only change the linear
combination by a multiple of $-1$, we define $\cV_n(M)$ to be the
subspace of $\C\cL_M$ spanned by all decorated links with $n$ double
points.  Clearly there is a filtration
\[\C\cL_M=\cV_0(M)\supset\cV_1(M)\supset\cdots\cV_n(M)\supset\cdots.\]
A map $f:\cL_M\rightarrow\C$ is a {\bf finite type invariant} if, when
extended linearly to $\C\cL_M$, $ \cV_{n+1}(M)$ lies in its kernel for
some $n$.  The {\bf order} of $f$ is the smallest such $n$.
  
Let $\iota : \cV_0 \rightarrow K(M)$ be the canonical projection.
\begin{prop}
$\iota(\cV_n) \subset h^nK(M)$.
\end{prop}

\begin{proof}
This is clearly true for $n=0$.  If 
$n>0$  we have
\begin{equation}
\begin{split}
\pm \dbp &= \rcr - \lcr \\
     &= -t \ift -t^{-1} \zer + t \zer +t^{-1} \ift \\
     &= (t-t^{-1}) \zer - (t-t^{-1})\ift,
\end{split}
\end{equation}
which implies $\iota(V_n) \subset h\iota(V_{n-1})$.  The result
follows by induction.  
\end{proof}

\begin{theorem}
If $K(M)$ is topologically free then the coefficients of a link are
finite type invariants.  Specifically, let $\Phi : K(M) \rightarrow
V[[h]]$ be an isomorphism; let $\{b_\alpha\;|\;\alpha\in\cI\}$ be a
basis for $V$; and define $\Phi_{i,\alpha}: K(M) \rightarrow \C$ by
\[ \Phi = \sum_{i}\left(\sum_\alpha \Phi_{i,\alpha} b_\alpha\right)h^i.\]
For each $i$ and each $\alpha$, $\cV_{i+1}(M)$ lies in the
kernel of $\Phi_{i,\alpha}\circ\iota$.
\end{theorem}

\begin{proof}
If $\beta \in \cV_{i+1}$ then $\Phi\circ\iota(\beta)\in h^{i+1}K(M)$.
Hence $\Phi_{i,\alpha}(\iota(\beta))=0$ for any $\alpha$. 
\end{proof}

\begin{rem} There are several different, though equally reasonable,
definitions of finite type invariants of unoriented, framed links.  In
each case, Proposition 4 and Theorem 2 hold.
\end{rem}

\subsection{Topological Generators}

As an algebra $K(F \times I)$ is generated by a finite set of knots
$\{K_1,\ldots, K_n\}$ \cite{B4}.  Ergo, as a module, $K(F \times I)$
is spanned by the set of all products of those generators.  Although
there is no product structure on $K(M)$ in general, it is possible to
capture the geometric essence of these products via the notion of
cabling.  For a framed knot $K$ in $M$ we define $K^n$ to be the
framed link formed by taking $n$ parallel copies of $K$, where
parallel means each copy is a pushoff of $K$ along its framing.  If
$n=0$ then $K^n=\emptyset$. If $L = K_1 \cup \cdots \cup K_m$ is a
framed link, let $L^{(n_1,\ldots,n_m)}$ be $K_1^{n_1} \cup \cdots \cup
K_m^{n_m}$.  We will refer to this as the $(n_1,\ldots,n_m)$-cable of
$L$.

For $L \in \cL_M$ let $W_L$ denote the $\C[[h]]$-linear span of all
cables of $L$ in $K(M)$.  If $W_L$ is dense in $K(M)$ we say that $L$
{\bf carries topological generators for $K(M)$}. The components of $L$
are the actual generators. A familiar example is $\emptyset$, carrying
generators for $K(S^3)$.  A non-example is the standard basis of
$K(F\times I)$, which consists of framed links that admit projections
with no crossings and no trivial components \cite{Pr}.  Unless $F$ is
planar with three or fewer boundary components, this basis contains
infinitely many knot types. Consequently, there is no link whose cables
span it.  

\begin{theorem}
Let $K_1, \ldots , K_n$ be knots that generate
$K_0(M)$ as an algebra. If $L_0$ is a link whose components are $K_1,
\ldots , K_n$, then $L_0$  carries topological generators for $K(M)$.
\end{theorem}

\begin{proof}
Consider the inductive hypothesis;
\[ H_n: \forall L \in K(M) \ \exists \beta \in W_{L_0} \text{\ such
that\ } \beta - L \in h^nK(M).\] Assume $H_n$ for all $i< n$, and
choose $L$. As a first step find $\beta_{n-1}$ so that $\beta_{n-1}-L
\in h^{n-1}K(M)$.  Let $\alpha$ be the constant term of a power series
representing $(\beta_{n-1}-L)/h^{n-1}$.  The idea is to write $\alpha$
as $\sum a_i L_i$, and resolve each link $L_i$ into cables of
$L_0$. To do this, recall that the knots $\{K_j\}$ generate
$K_0(M)$. That means there exists $\gamma_i \in W_{L_0}$ such that
$\gamma_i - L_i \in hK(M)$. It follows that there exists $\beta_n \in
W_{L_0}$, $\beta_n +\beta_{n-1} - L \in h^nK(M)$. Hence $H_n$ is true
for all $n$. Since $\{h^nK(M)\}$ is a neighborhood basis for $0$, we
have shown that $L$ lies in the closure of $W_{L_0}$. 
\end{proof}

\begin{sch}
If $K_0(M)$ is finitely generated as a $\C$-module, then $K(M)$ is finitely generated as a $\C[[h]]$-module.
\end{sch}

\begin{proof}
Let $\alpha_1,\ldots,\alpha_n$ be elements of $K(M)$ such that $K_0(M)$
is generated by their images under the canonical projection.  Modify
the inductive hypothesis of the preceding proof to read:
\[ H_n: \forall L \in K(M) \ \exists \beta \in \text{\ span}(\{\alpha_i\} 
\text{\ such that\ } \beta - L \in h^nK(M).\] The proof that $H_n$
holds for all $n$ is essentially the same as above.
\end{proof}

\begin{cor}  If $K_0(M)$ has no nilpotent elements then $X[M]$ induces topological generators for $K(M)$.  More specifically, there exists a finite set $tr_{\gamma_1},\ldots,tr_{\gamma_n}$ which generates $X[M]$ as a ring.  If $K_1,\ldots,K_n$ are knots chosen to correspond to each $tr_{\gamma_i}$, then they are topological generators of $K(M)$.
\end{cor}

Recall that a 3-manifold is {\bf small} if it contains no
incompressible surface.

\begin{theorem}
If $M$ is a small 3-manifold, then $K(M)$ is finitely generated as a
module.
\end{theorem}

\begin{proof}
By Theorem 11 of \cite{B1}, $K_0(M)$ is a finitely generated module.
\end{proof}

It is interesting to compare Theorem 4 to a conjecture of Traczyk
 \cite[Problem 1.92 H]{Ki}.

\begin{conjecture}  If $M$ is a simply connected $3$-manifold
other than $S^3$ (possibly with holes), then $S_{2,\infty}(M)$
is infinitely generated.
\end{conjecture}

 The module $S_{2,\infty}(M)$ is
the Kauffman bracket skein module defined over
${\mathbb Z}[A,A^{-1}]$.  For any $3$-manifold $M$,
there is a map $S_{2,\infty}(M) \rightarrow K(M)$ obtained
by substituting $-t$ for $A$.  Is this map injective?
If it isn't injective what does its kernel look like?

\end{document}